\documentclass[letterpaper, 10 pt, conference]{ieeeconf}  

\IEEEoverridecommandlockouts
\usepackage{cite}
\usepackage{mdframed}
\usepackage{soul}
\usepackage[inkscapelatex=false]{svg}
\usepackage{bm}

\makeatletter
\newcommand{\subalign}[1]{%
  \vcenter{%
    \Let@ \restore@math@cr \default@tag
    \baselineskip\fontdimen10 \scriptfont\tw@
    \advance\baselineskip\fontdimen12 \scriptfont\tw@
    \lineskip\thr@@\fontdimen8 \scriptfont\thr@@
    \lineskiplimit\lineskip
    \ialign{\hfil$\m@th\scriptstyle##$&$\m@th\scriptstyle{}##$\hfil\crcr
      #1\crcr
    }%
  }%
}
\makeatother

\newmdenv[linecolor=green, backgroundcolor=green!10]{hlgb}
\newmdenv[linecolor=red, backgroundcolor=red!10]{hlrb}

\def\BibTeX{{\rm B\kern-.05em{\sc i\kern-.025em b}\kern-.08em
    T\kern-.1667em\lower.7ex\hbox{E}\kern-.125emX}}
\usepackage{xcolor}

\usepackage[utf8]{inputenc}
\usepackage[english]{babel}
\usepackage[normalem]{ulem}

\usepackage{mathptmx} 
\usepackage{amsmath} 
\usepackage{amssymb}
\usepackage{amsfonts}
\usepackage[]{algorithm2e}
\SetKwInput{kwInput}{Input}
\usepackage{esvect}

\usepackage{mathtools}

\usepackage{bm}
\usepackage{cite}
\usepackage{comment}
\allowdisplaybreaks

\numberwithin{theorem}{section}

\usepackage{graphicx}
\usepackage{subfig}
\usepackage{epstopdf}
\usepackage{hyperref} 
\usepackage{color}
\definecolor{Questions}{HTML}{1F77B4}

\usepackage{cancel}

\usepackage{tikz}
\usetikzlibrary{calc}
\usetikzlibrary{arrows,positioning} 
\usetikzlibrary{intersections, patterns}
\usetikzlibrary{shapes}
\usepackage{mathtools}
\usepackage[cmtip,all]{xy}
\newcommand{\longsquiggly}{\xymatrix{{}\ar@{~>}[r]&{}}}

\DeclareMathAlphabet{\mathcal}{OMS}{cmsy}{m}{n}

\begin{document}
\title{Navigation with shadow prices to optimize multi-commodity flow rates}

\author{Ignacio Boero$^{1, 2}$, Igor Spasojevic$^{1}$, Mariana del Castillo$^{2}$ George Pappas$^{1}$, Vijay Kumar$^{1}$, Alejandro Ribeiro$^{1}$ 
\thanks{This work was supported in part by NSF Grant CCR-2112665. $^{1}$Authors are with the Department for Electrical and Systems Engineering,
        University of Pennsylvania, USA; 
        $^{2}$Authors with the Institute of Electrical Engineering, Universidad de la República, Uruguay; 
        {\tt\small iboero, mdelcastillo@fing.edu.uy, igorspas, pappasg, kumar, aribeiro@seas.upenn.edu}}%
}

\maketitle

\begin{abstract}
We propose a method for providing communication network infrastructure in autonomous multi-agent teams. 
In particular, we consider a set of communication agents that are placed alongside regular agents from the system in order to improve the rate of information transfer between the latter.
In order to find the optimal positions to place such agents, we define a flexible performance function that adapts to network requirements for different systems. 
We provide an algorithm based on shadow prices of a related convex optimization problem in order to drive the configuration of the complete system towards a local maximum. 
We apply our method to three different performance functions associated with three practical scenarios in which we show both the performance of the algorithm and the flexibility it allows for optimizing different network requirements.
\end{abstract}
\section{Introduction}

Autonomous multi-agent systems have lately found applications in numerous challenging tasks. 
Examples include a team of robots mapping an unknown or disaster-stricken environment, performing surveillance missions, and carrying out search and rescue operations. 
A key hurdle lies in the fact that individual agents have to make decisions with outdated or even missing information from their teammates. 
This is due to the inherent problem that communications cannot be instant. 
To overcome this issue, a large body of research has been done on developing algorithms to take optimal decisions with the available information. 
A parallel line of works seeks to improve the communication network infrastructure formed by a multi agent system, with the aim of making communications faster and more reliable. 

This paper belongs to the second line of research. 
A common approach to solving this problem involves deploying a second team of agents whose purpose is to provide network infrastructure for the first team \cite{mox2020mobile,4543418, stephan2017concurrent, fink2013robust,6213085}. 
The problem then amounts to positioning the newly deployed agents in order to allow the greatest amount of communication flow. 
Previous work has predominantly focused on maximizing the connectivity of the graph defined by the spatial configuration of agents \cite{zavlanos2011graph, kim2005maximizing, zavlanos2005controlling,mox2022learning}. 
This line of works relies heavily on quantities from algebraic graph theory that capture a notion of connectivity. 
Indeed, one common figure of merit for the connectivity of a graph is the second eigenvalue of its Laplacian.
The latter can be computed efficiently, and many approaches solve the problem by finding methods to maximize the second eigenvalue. 
Although such a notion of connectivity is attractive in that it captures a complex notion with a single number, we believe that resulting solutions often suffer from lack of flexibility to adapt to different communication requirements that a multi agent system needs. 
In our work, we consider a solution that captures a more complete notion than only of connectivity. This has been previously done by \cite{6213085}, nonetheless we differ in that we depart from the assumption that we know the minimum desired communication rate for each agent, which many times is not the case. This modifies the framework from restricting agents movement to feasible configurations to moving in directions that maximize an objective.

In particular, we define the figure of merit for a given configuration of agents using the solution of a Multi-Commodity Flow Problem. 
Multi-Commodity Flow Problems (MCFP) are widely used in many areas \cite{feo1989flight,foulds1981multi}, including network modeling \cite{assad1978multicommodity,ali1984multicommodity}. 
It consists of a graph and a set of commodities, where a commodity can model the ``supply'' of information from some nodes, usually named source nodes, to other ``demand'' nodes, named sink nodes. 
The restrictions on this problem are that there is a cost associated for a commodity to flow from one node to another, and edges between nodes have limitations on how much cost they can withstand. 
Then, the problem is to maximize a utility function that increases with the amount of commodities that can be pushed from sources to sinks, while respecting given edge capacities. 
This model bears a straightforward translation to our problem. 
In our setting, the graph is defined by the spatial configuration of the agents, the commodities are information being transmitted between agents, and the capacity constraints between nodes are the link capacities between agents. 
Then, the performance of a configuration is given by the value of the utility function for the solution of that MCFP. 
The flexibility of this figure of merit lies in the fact that different agents can be selected as source or sink agents.
Furthermore, the options for the utility function can capture a range of different modelling aspects that are more refined than mere graph connectivity. 

Having defined the figure of merit, we present an algorithm that from random initial positions moves select agents to positions that locally maximize the utility of the MCFP. 
When moving an agent, we perturb the capacity constraints of the MCFP, as we increase the link between the moved agent and the agents to whom the distance decreased, and vice-versa. 
Then, the dilemma is to choose which link to increase in order to maximize the value of the MCFP. 
For doing that we rely on perturbation theory, which establishes that dual variables associated with a constraint are closely related to the gain of relaxing that constraint \cite{boyd2004convex}. 
We then do numerical experiments which demonstrate the efficacy of the algorithm, as it consistently increases the utility function of the MCFP on all experiments.
Ultimately, we show the flexibility of this approach, illustrating how different practical scenarios can lead to different definitions for the MCFP, which in turn leads to different solutions of the algorithm. 

Regarding the rest of the paper, in section II we give the mathematical formulation for a general figure of merit defined by a MCFP. 
In section III we explain the proposed method, and provide a mathematical justification. 
In section IV we apply the algorithm to different use cases involving various practically-motivated static scenarios, 
and we finish with a test of our approach on a dynamic scenario. 
Finally, in section V we summarize our findings and propose future work.

\section{Dynamic Multi-commodity Flow Problem}

We consider a maximum multi-commodity flow problem on a dynamic graph $\mathcal{G} = (\mathcal{V}, \mathcal{E})$, with vertex set $\mathcal{V}$ and edge set $\mathcal{E}$. 
Vertices in $\mathcal{V}$ correspond to agents in a robot team whose configurations determine capacities of edges in $\mathcal{E}$.
Letting $|\mathcal{V}| = N$, we identify $\mathcal{V}$ with $[N] = \{1,2,...,N\}$, and denote the configuration of agent $i \in [N]$ by $\mathbf{x}_i \in \mathbb{R}^2$.
The capacity function
\begin{equation}
\begin{aligned}
c : \mathbb{R}^2 \times \mathbb{R}^2 & \rightarrow \mathbb{R}_{++} \\
\end{aligned}
\end{equation}
constrains the maximum amount of cumulative flow along any edge $(i,j) \in \mathcal{E} \subseteq \mathcal{V} \times \mathcal{V}$ by $c(\mathbf{x}_i, \mathbf{x}_j)$.
Here, $c$ is a positive differentiable function, with $c(\mathbf{x},\mathbf{y})$ decreasing in $|| \mathbf{x} - \mathbf{y} ||_2$ for fixed $\mathbf{x}$ (or $\mathbf{y}$) and fixed $(\mathbf{x}-\mathbf{y})/||\mathbf{x}-\mathbf{y}||_2$.
Letting $\mathcal{I} \subseteq \mathcal{V}$ denote agents with controllable configurations,  we solve: 
\begin{equation} \tag{DMCF}
\max_{(\mathbf{x}_i)_{i \in \mathcal{I}}} \Phi(\mathbf{x}_{1:N})
\label{eqn:global_problem}
\end{equation}
where $\Phi$ represents the value of the following \textit{generalized} multi-commodity flow problem 

\begin{equation}
\tag{P-MCF}
\begin{aligned}
&\Phi(\mathbf{x}_{1:N}) \  =  \max_{\textbf{r} \in[0,1]^{N \times N \times K}, \ \mathbf{a} \in \mathbb{R}_+^{N \times K}} \quad \mathcal{U}\left((a_i^k)_{\ k \in \mathcal{A}, \ i \in S_k} \right)\\
s.t. & \\
& a^k_i \leq   \sum_{j=1}^{N} r_{i,j}^{(k)}  -  \sum_{j=1}^{N} r_{j,i}^{(k)} \quad \forall k \in \mathcal{A}, \ \forall  i \in \mathcal{S}_k  \\
& \sum_{j=1}^{N} r_{i,j}^{(k)}  -  \sum_{j=1}^{N} r_{j,i}^{(k)}  = 0 \quad \forall k \in \mathcal{A}, \ i \in \mathcal{I} \\
&\sum_{k \in \mathcal{A}} r_{i,j}^{(k)} \leq c(\mathbf{x}_i, \mathbf{x}_j) \quad \forall i,j \in \mathcal{V}. \\
\end{aligned}
\label{eqn:P-MCF}
\end{equation}
In the latter, $\mathcal{A} = \mathcal{V} \setminus \mathcal{I}$ represents the subset of agents that need to exchange information. 
We have one commodity for each agent in $\mathcal{A}$, indexed by $k \in \mathcal{A}$, and let $K := |\mathcal{A}|$ denote the total number of commodities. 
We also let $I := |\mathcal{I}|$, and note that $N = K + I$.
For commodity $k$, agents in $S_k \subseteq \mathcal{A}$ act as source nodes while the rest of the agents in $\mathcal{A}$ act as sinks. 
This is (implicitly) encoded by the first constraint of Problem \ref{eqn:P-MCF}.
In any flow, agents in $\mathcal{I}$ are relay nodes.
They neither generate nor accumulate information, as described by the second constraint of the problem. 
The edge capacity between any two agents represents the maximum rate at which they can send information between themselves across all flows. 
The latter is captured by the third constraint of Problem \ref{eqn:P-MCF}.
Ultimately, the objective of \ref{eqn:P-MCF} is a general one in that we only require concavity of the global \textit{utility} function $\mathcal{U}$.
Classical maximum multi-commodity flow problems focus on cases where $\mathcal{U}$ is the sum of  $(a^{k}_i)_{k \in \mathcal{A}, i \in \mathcal{S}_k}$. 
The present formulation allows us to generalize this to other measures of team performance, as we will see in section IV.
Our problem is dynamic due to edge capacities being a function of the positions of the agents, which can be varied. 
This differs from usual multi-commodity flow problems in which the capacity constraints are fixed, and only the flows can be optimized.
The set of decision variables therefore consists of positions of agents in $\mathcal{I}$, denoted by $(\mathbf{x}_i)_{i \in \mathcal{I}}$, as well as  $(r_{i,j}^{(k)})_{i,j \in \mathcal{V}, k \in \mathcal{A}}$, the flow for each commodity $k$ being transmitted between each pair of agents $(i,j)$, and $(a^{k}_i)_{k \in \mathcal{A}, i \in \mathcal{S}_k}$, the amount of commodity $k$ generated by agent $i$.
We also introduce vector notation for the latter two variables, where $\mathbf{r} \in [0,1]^{N\times K\times K}$ collects $r_{ij}^{(k)}$ in ascending lexicographic order, first by source $i$, then by sink $j$, and finally by commodity $k$, while $\mathbf{a} \in \mathbb{R}_+^{N\times K}$ arranges $a^{k}_i$ similarly, first by agent $i$ and then by commodity $k$.

\section{Shadow Price Ascent Algorithm}

\subsection{Proposed Algorithm}
As $c(\mathbf{x},\mathbf{y})$ cannot be a concave function, \eqref{eqn:global_problem} is in general a non-convex problem. 
As a result, we develop an algorithm for finding a first order stationary point.
We use an iterative method where at each step we solve the dual problem of \ref{eqn:P-MCF}, and update the positions of $(\mathbf{x}_{i})_{i \in \mathcal{I}}$ along a direction of local increase.
The dual of \ref{eqn:P-MCF} is given by:

\begin{equation}
\tag{D-MCF}
\begin{aligned}
D^* (\mathbf{x}_{1:N}) &=\min_{\subalign{&\bm{\lambda} \in \mathbb{R}_+^{|S_k|\times K}, \\ &\bm{\mu} \in \mathbb{R}_+^{K\times I} , \\ &\bm{\nu} \in \mathbb{R}^{N \times N}}} D(\bm{\lambda},\bm{\mu},\bm{\nu}, \mathbf{x}_{1:N}) \\
&=\min_{\subalign{&\bm{\lambda} \in \mathbb{R}_+^{|S_k|\times K}, \\ &\bm{\mu} \in \mathbb{R}_+^{K\times I} , \\ &\bm{\nu} \in \mathbb{R}^{N \times N}}}  \ \max_{\subalign{&\textbf{r} \in[0,1]^{N \times N \times K} \\ &\mathbf{a} \in \mathbb{R}_+^{N \times K}}}  \ L(\bm{\lambda},\bm{\mu},\bm{\nu}, \mathbf{a},\mathbf{r},\mathbf{x}_{1:N}),
\end{aligned}
\label{eqn:D-MCF}
\end{equation}
where $L$ is the Lagrangian of \eqref{eqn:P-MCF} defined as:
\begin{equation}
\begin{aligned} 
L(\bm{\lambda},\bm{\mu},\bm{\nu},\bm{a},\bm{r},\mathbf{x}_{1:N}) = & \quad  \mathcal{U}\left((a_i^k)_{\ \forall k \in \mathcal{A}, \ i \in S_k} \right) \\ 
&-\sum_{i\in \mathcal{S}_k,k \in \mathcal{A}} \lambda^{k}_{i} [a_{i}^{k} - (\sum_{j=1}^{N} r_{i,j}^{k}  -  \sum_{j=1}^{N} r_{j,i}^{k})]\\
&-\sum_{i \in \mathcal{I},k \in \mathcal{A}} \nu^{k}_{i} [\sum_{j=1}^{N} r_{i,j}^{(k)}  -  \sum_{j=1}^{N} r_{j,i}^{(k)} ]\\
&-\sum_{i,j \in \mathcal{V}} \mu_{i,j}[(\sum_{k \in \mathcal{A}} r_{i,j}^{(k)}) - c(\mathbf{x}_i, \mathbf{x}_j)].
\end{aligned}
\label{eqn:L-MCMF}
\end{equation}
The dual variables $ \bm{\lambda} \in \mathbb{R}_+^{|S_k|\times K}, \ \bm{\nu} \in \mathbb{R}^{N \times N}, \ \bm{\mu} \in \mathbb{R}_+^{K\times I}$ are associated with constraints one, two, and three, respectively.

In light of the fact that $c > 0$, it is not difficult to show that  \eqref{eqn:P-MCF} is a convex program for which Slater's condition holds. 
In particular, it has zero duality gap. 
Solving the dual problem \ref{eqn:D-MCF} is therefore equivalent to solving the primal. 
Nevertheless, for our purposes it is more convenient as it allows us to readily extract the sensitivity of the objective of the team problem \ref{eqn:global_problem} with respect to configurations of controllable agents. \\

\textbf{Proposition 1.} Let $\{\mu_{ij}^*(x_{1:N})\}_{i,j \in \mathcal{V}}$ be the solution to \eqref{eqn:D-MCF} for a given set of positions $\mathbf{x}_{1:N}$.
Then $ \partial \Phi / \partial X \in (\mathbb{R}^2)^{\mathcal{I}} $ given by
\begin{equation}
\left. \left( \frac{\partial \Phi}{\partial X} \right)_i  \right|_{\mathbf{x}_{1:N}} = \sum_{j \in \mathcal{V}} (\mu_{ij}^*(\mathbf{x}_{1:N}) + \mu_{ji}^*(\mathbf{x}_{1:N}))\nabla_{\mathbf{x}_i} c(\mathbf{x}_i,\mathbf{x}_j)|_{\mathbf{x}_{1:N}}
\label{eqn:subG-P}
\end{equation}
with $i \in \mathcal{I}$ is a direction of local increase for \ref{eqn:global_problem}.\\

\begin{algorithm} 
\RestyleAlgo{ruled}
\caption{Shadow Price Ascent}  \label{alg:alg2}
\kwInput{$X_{\mathcal{V} \setminus \mathcal{I}}$, $X_{\mathcal{I}}^{0}$, $\alpha$, $tol$}
$X_{\mathcal{I}} \gets X_{\mathcal{I}}^{0}$\\
$\Phi_{0} \gets -\infty, \ \Delta \Phi \gets +\infty$\\
$(C)_{ij} \gets c(x_i,x_j) \ \forall x_i,x_j \in X_{\mathcal{V}}$\\
$(\partial C)_{ij} \gets \nabla_{x_i} c(x_i,x_j) \ \forall x_i,x_j \in X_{\mathcal{V}}$\\
$t \gets 0$\\
\While{$\vert \Delta \Phi \vert \geq$ tol}{
 $(\mu^{*}, \Phi_{t+1}) \gets$ Solve-DMCF$(C)$\\
 \For{$i \in \mathcal{I}$}{
    $\left( \frac{\partial \Phi}{\partial X} \right)_{i} \gets  \sum_{j \in \mathcal{V}} (\mu_{ij}^* + \mu_{ji}^*)(\partial C)_{ij}$\\
 }
 $X_{\mathcal{I}} \gets X_{\mathcal{I}} + \alpha* \left( \frac{\partial \Phi}{\partial X} \right)_{\mathcal{I}}$\\
 $(C)_{ij} \gets c(x_i,x_j) \ \forall x_i,x_j \in X_{\mathcal{V}}$\\
 $(\partial C)_{ij} \gets \nabla_{x_i} c(x_i,x_j) \ \forall x_i,x_j \in X_{\mathcal{V}}$\\
 $t \gets t+1$ \\
 $\Delta \Phi \gets \Phi_{t+1} - \Phi_{t}$\\
 }
\end{algorithm}

Our approach is summarized in Algorithm \ref{alg:alg2}. 
From a high level, it is a gradient ascent method. 
Lines $1-5$ initialize all variables.  
Thereafter, until the change in the team objective function $\Phi$ drops below the specified tolerance value $tol$, the algorithm calculates a local direction of ascent by using the dual variables in an optimal solution to problem \ref{eqn:D-MCF} obtained by the interior point solver MOSEK \cite{mosek}.  
We refer to the latter sub-procedure by the method ``Solve-DMCF''; it takes as input a graph with fixed edge capacities, and outputs the optimal flows along every edge, the team value function, as well as the optimal dual (Lagrangian) variables for every constraint. 
Lines $10-12$ then leverage the dual variables to perform the gradient ascent step using a stepsize $\alpha$ that is tuned for the family of problems at hand.

\subsection{Intuition behind the direction of local increase}

Any algorithm that seeks to maximize the cumulative capacity of the network connecting the agents faces a trade-off.
In order to maximize the rate of information transfer between agent $i$ and $j$, it needs to update the position of agent $i$ in the direction  of the gradient of $c(\mathbf{x}_i,\mathbf{x}_j)$ w.r.t $\mathbf{x}_i$.
Nonetheless, it is easy to see that this heuristic can be problematic, as maximizing the quality of the connection between one pair of agents could easily diminish the quality between another. 
A simplistic approach would be to move agent $i$ along the average of the gradients of $c(\mathbf{x}_i,\mathbf{x}_j)$ w.r.t $\mathbf{x}_i$ across all $j$.
However, this approach would miss important information regarding the complexity of how information flows through the network. 
In particular, if a link is a bottleneck in the network, it would be desirable to increase its capacity  before that of others.
Our algorithm attacks this problem via leveraging a well known concept on constrained optimization - the sensitivity of the perturbation function.
In a constrained optimization problem in which strong duality holds, the dual variable associated with a constraint can be seen as the sensitivity of the primal function to perturbations of the value of the bound in that constraint.
Bottlenecks manifest as the edges to which the primal function displays maximal sensitivity, thereby emerging as most influential during the update of agent positions.

\section{Numerical Experiments}

\begin{figure}[h]
    \centering
    \includegraphics[scale=0.25]{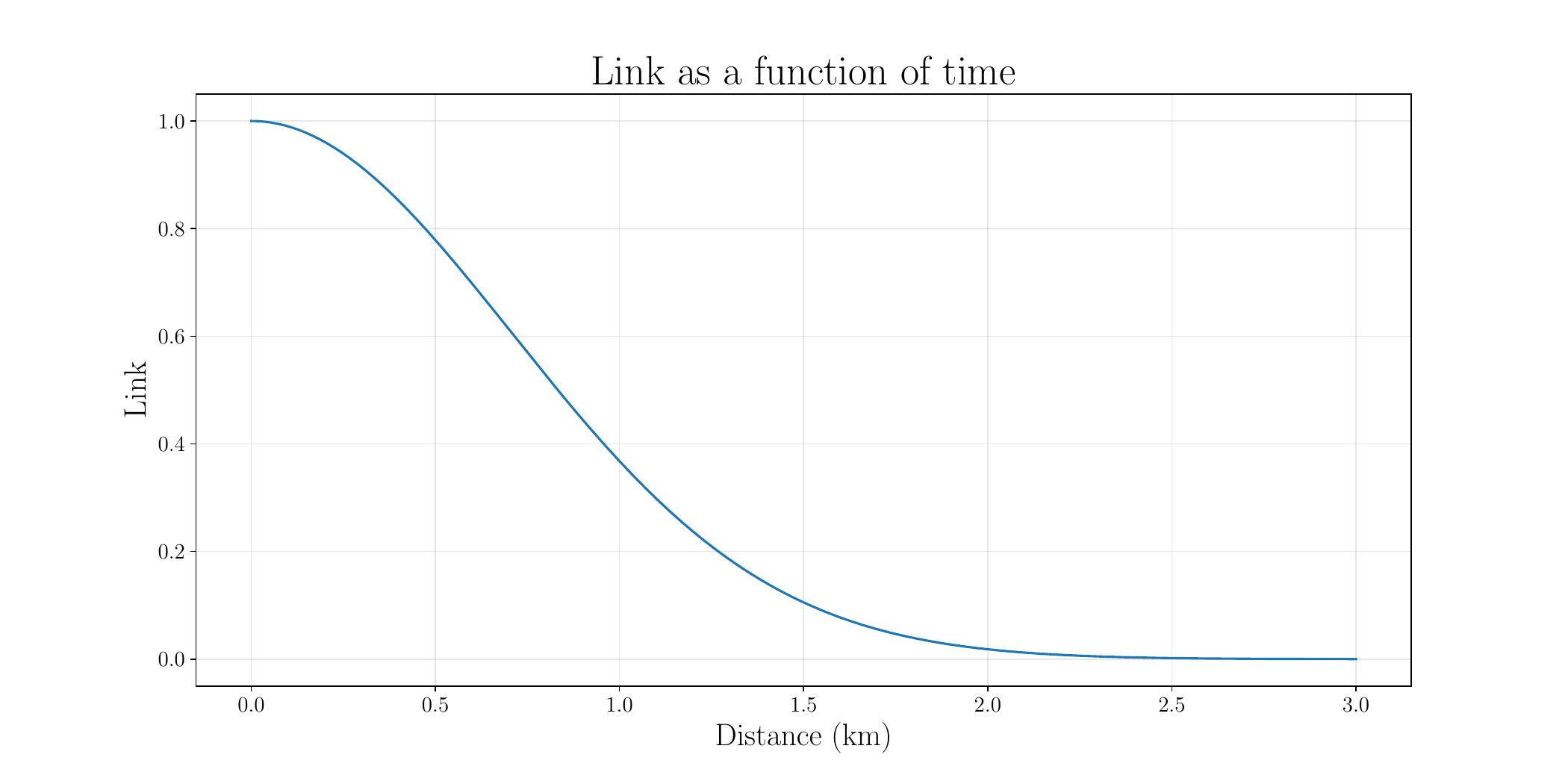}
    \caption{Function $c(\mathbf{x},\mathbf{x})$ from  equation \eqref{eqn:c_func} with $D = 2$ and $d_0$ = 1km}
    \label{fig: cfunc}
\end{figure}

In this section, we provide numerical experiments to show the performance of our algorithm. 
First, we set several parameters that will be used throughout all experiments.
We use the capacity function 
\begin{equation}
c(\mathbf{x}, \mathbf{y}) = exp \left( -  \left( \frac{|| \mathbf{x} - \mathbf{y} ||_2}{d_0} \right)^{D}  \right)
\label{eqn:c_func}
\end{equation}
Intuitively, $D$ captures the shape of the fading of the capacity of the link between a pair of agents as function of their distance, whereas $d_0$ captures the characteristic length scale of the decay.
In all experiments, the density of agents $\rho$ is the same; the size of the area in which agents are spawned will depend on the number of agents. 
We set $\rho = 1 \ agents/km^2$, $D=2$, and $d_0 = 1km$.  
A plot of the rate function is shown in Figure \ref{fig: cfunc}.
We used the step size $\alpha = 0.4$, scaling it by the value of $0.97$ after every iteration.
For all experiments we restrict the utility function of the P-MCF to a linear combination of the minimum $a^k_i$ of each commodity $k$. 
Therefore it is defined as 
$$ \mathcal{U}\left((a_i^k)_{\ k \in \mathcal{A}, \ i \in S_k} \right) = \sum_{k \in \mathcal{A}} \omega_k \min_{i \in S_{k}} a^k_i$$

We present three distinct study cases where the merit function above can be used to model different network requirements.
\subsubsection{Ad-Hoc Networking}
The agents represent a mobile Ad-Hoc Network; they navigate together and need to remain connected with as much network capacity as possible. 
In this case, we set $\omega = \mathbf{1}$, weighing all commodities equally.
\subsubsection{Routing to infrastructure AP}
A set of agents is performing a task where they need to maintain communication with an access point (AP). 
This setting can arise in scenarios where agents carry out a search and rescue mission in a remote area, and need to relay local information to an AP. 
There is only one sink, and so only one commodity is required. 
Therefore, $\omega$ is set to be all zero except for the $j$-th entry, corresponding to the AP, which is set to 1. 
\subsubsection{Distributed Algorithm}
Agents are running a distributed algorithm that periodically needs to update information from surrounding agents. 
We set $\omega$ to be a time-varying vector, whose entries associated with agents that need to update their information for a given instant are set to one, with remaining entries set to zero. 

In the remaining part of this section, we use the previous three study cases to show multiple properties of our method.

\subsection{Performance}
\begin{figure}[h]
    \centering
    \includegraphics[scale=0.25]{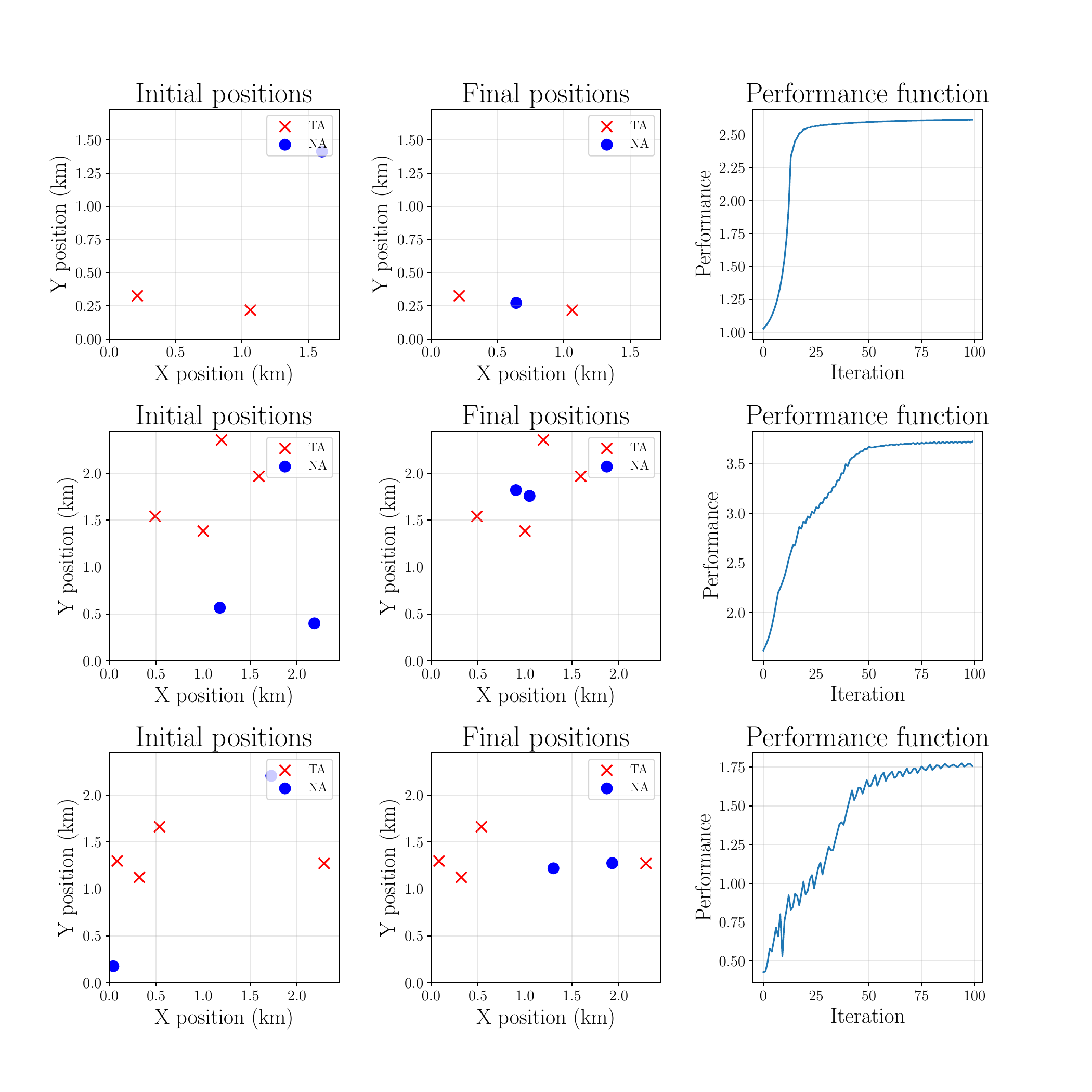}
    \caption{Test with small teams of task and network agents in various configurations. The first row involves two task agents and one network agent. The second row has four task agents arranged in a square and two network agents. The third row corresponds to a scenario with four task agents, where one is separated from the other three, and two network agents. In all cases, the first column depicts initial positions of all agents, the second column the positions of networks agents adjusted by the algorithm, and the third column the evolution of the team utility function with the number of iterations.}
    \label{fig:intuitive}
\end{figure}

We start by showing results for several scenarios with a small number of agents, in which the form of optimal solutions can be gauged by inspection. 
In all such setups, we use the first study case, setting $\omega = \mathbf{1}$.
We reserve the term ``network agents'' for agents in the controllable subset $\mathcal{I} \subseteq \mathcal{V}$, and the term ``task agents'' for those in subset $\mathcal{V} \setminus \mathcal{I} \equiv \mathcal{A}$.
First, we consider a problem with two task agents and one network agent. 
The expected result is for the network agent to be placed in the middle of the segment formed by its teammates. 
The next scenario is composed of four task agents forming a square shape, and two network agents. 
Here, the intuitive solution is to place the latter pair close to the middle of the square.
The last scenario also involves four task agents; now, however, three of them are close together, while the fourth one is far away. 
The last example serves to illustrate an important aspect of the utility function. 
We defined it to be the sum of $a^k$, where $a^k$ is the \textit{minimum} $a^k_i$ across all $i \in S_{k} \ (\text{in this case } S_{k} = \mathcal{A} \setminus \{k\})$ for commodity $k$. 
Therefore, if one agent is disconnected from the rest, it will lead to $a^k=0$. 
For this reason, we expect the algorithm will connect this agent to its teammates, in order to increase $a^k$. 
Figure \ref{fig:intuitive} illustrates outputs of the algorithm running in these three scenarios. 
In all cases, expected behaviour emerges. 
Furthermore, the team utility function increases with the number of iterations, thus providing reassuring evidence of the soundness of the algorithm.

\subsection{Flexibility}

\begin{figure}[h]
    \centering
    \includegraphics[scale=0.25]{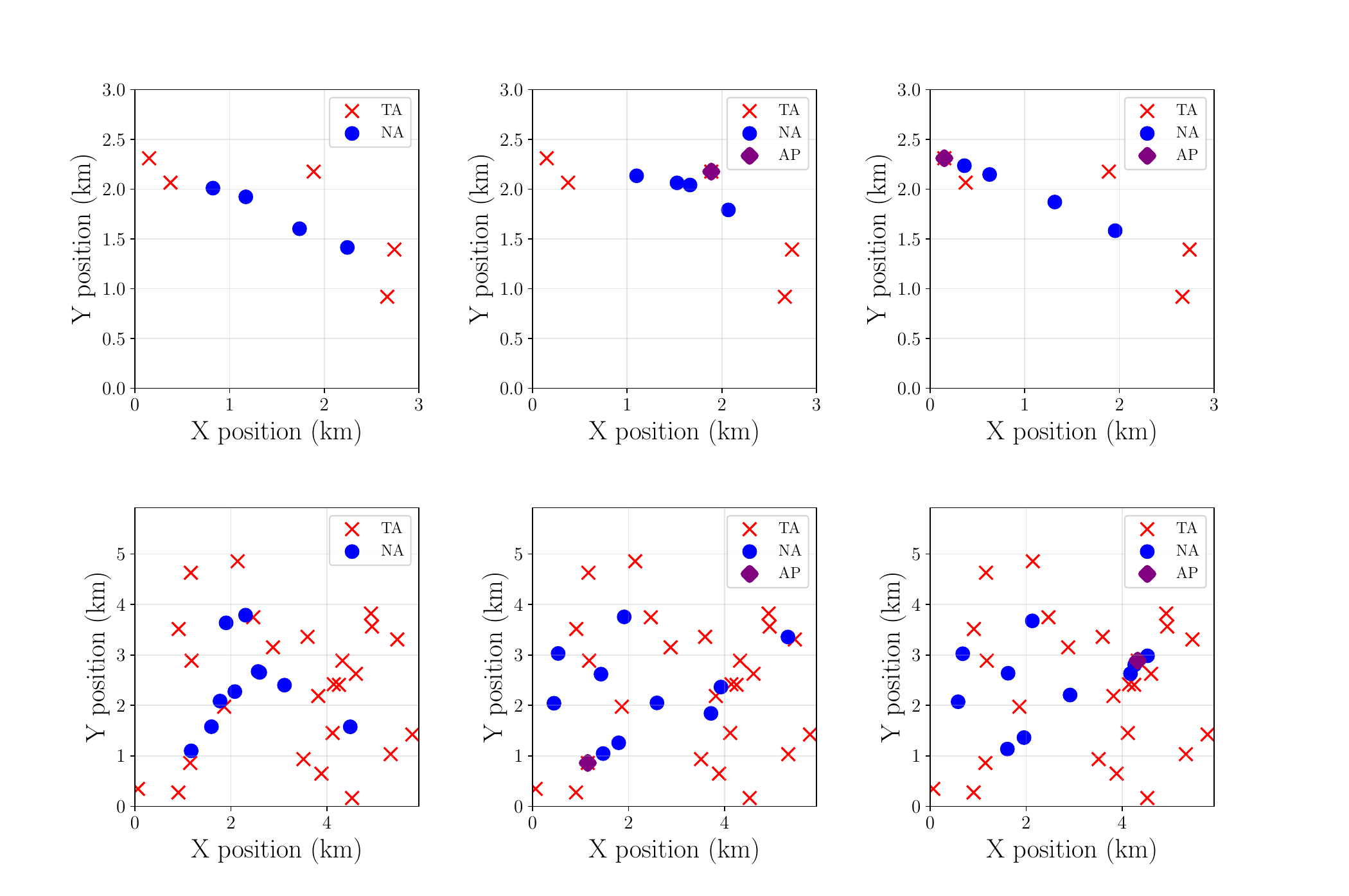}
    \caption{Test with one commodity in the objective. The first row corresponds to 5 task agents and 4 network agents. The second row involves 25 task agents and 10 network agents. In both cases, the first column is the solution with all weights equal to one. In the second and third column, the weight of the commodity of one randomly selected task agent was set to one, while the rest were set to zero.}
    \label{fig:case1vscase2}
\end{figure}

We consider two scenarios, in each of which we run the algorithm three times, using the same task agent configuration, but varying the utility function. 
In the first run of the algorithm, we use the all-ones weight vector, simulating the first study case. 
For the next two runs, simulating the second study case, we choose a different task agent to play the role of the AP. 
In the first scenario, we use 5 task agents and 4 network agents, and in the second scenario we use 25 task agents and 10 network agents.
The positions of network agents computed by our algorithm are shown in Figure \ref{fig:case1vscase2}. 
We can observe a marked difference in solutions for each case. 
While using $\omega = \mathbf{1}$, the solution exhibits higher spatial dispersion, in line with the intuitive goal to connect all agents. 
However, when there is only one commodity in the objective, the communication agents are placed closer to the agent playing the role of the AP.

\begin{figure}[h]
    \centering
    \includegraphics[scale=0.25]{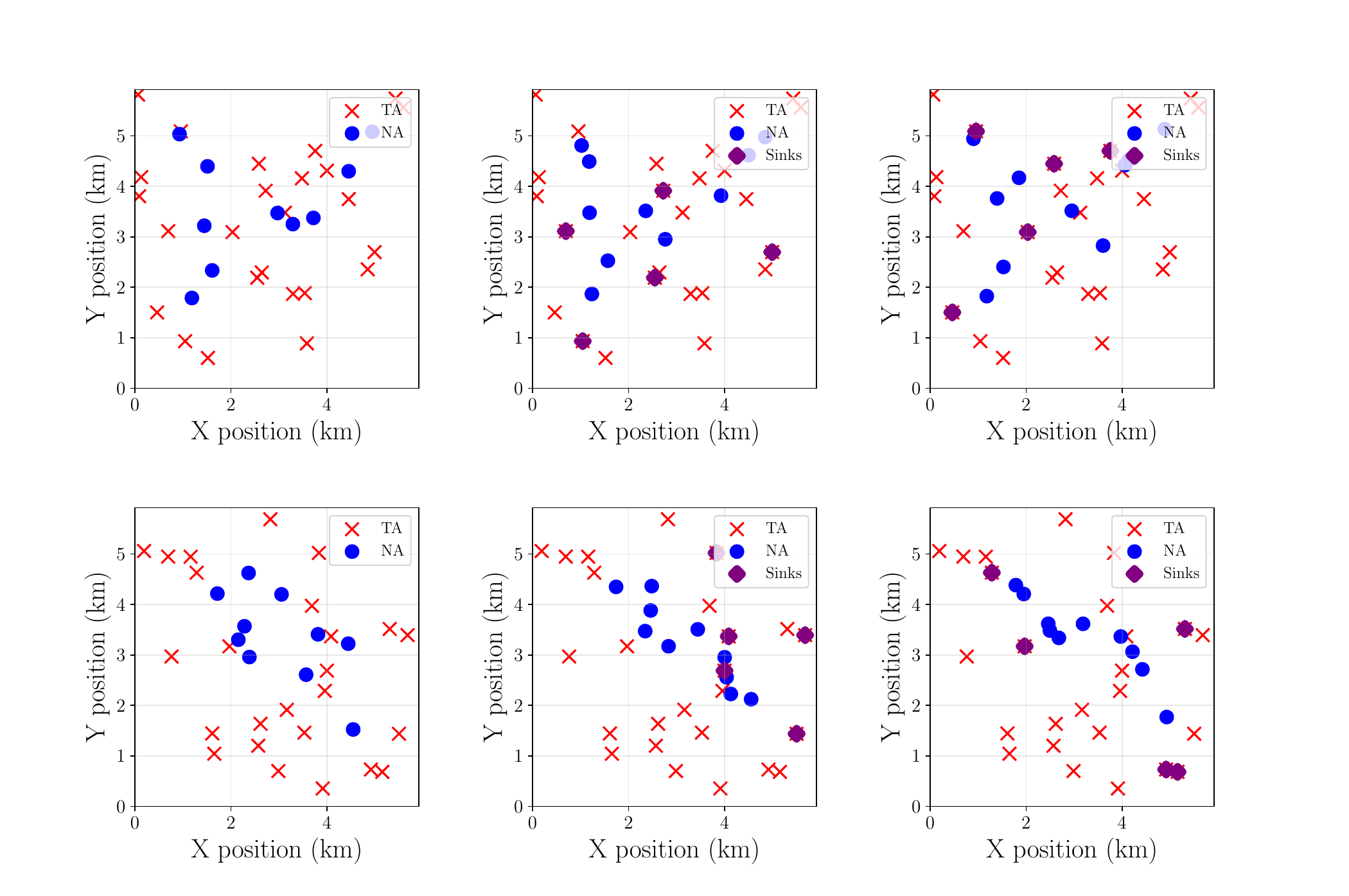}
    \caption{Test comparing different utility functions. Both rows illustrate examples with 25 task agents and 10 communication agents. In the first column all 25 commodities are active, whereas in the second and third 5 out of the 25 commodities are active.}
    \label{fig:mulitple-sinks}
\end{figure}

We also repeat the latter experiment, this time comparing the first case study with the third. 
Figure \ref{fig:mulitple-sinks} shows how, once again, different solutions are reached, where in cases where commodities are restricted to a few task agents, the communication agents end up placed close to them.

\subsection{Scalability}

\begin{figure}[h]
    \centering
    \includegraphics[scale=0.25]{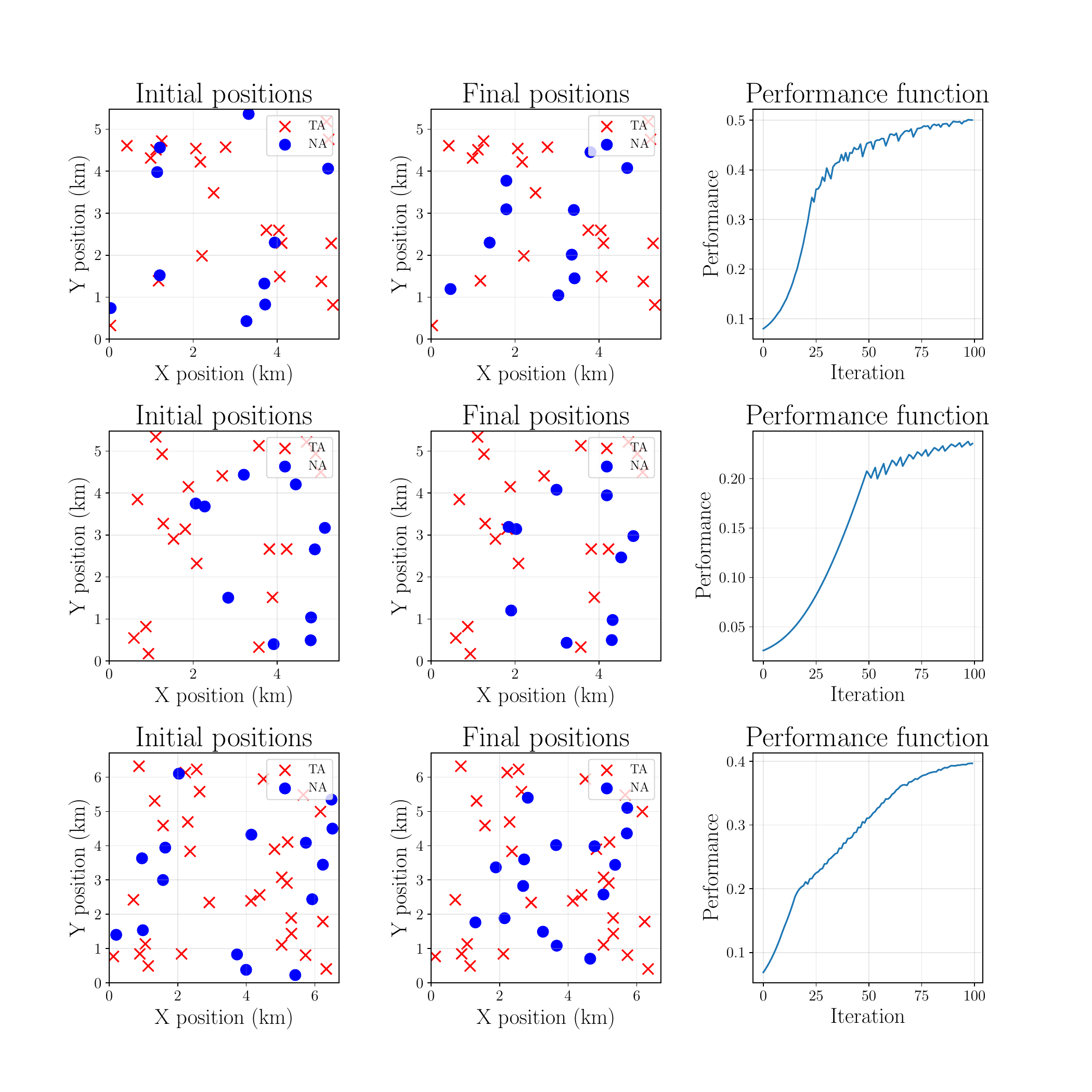}
    \caption{Test with weight vector equal to one. The first row involves 20 task agents and 10 network agents, and the second involves 25 task agents and 10 network agents. The third row corresponds to 30 task agents and 15 network agents. In all cases the first column depicts initial positions of all agents, the second column the output of the algorithm, and the third column the evolution of the performance function.}
    \label{fig:scalability}
\end{figure}

Our final goal was to determine how the computational performance of our method scales with the number of agents/size of the graph. 
To this end, we ran 100 simulations for different numbers of task agents, and computed the average time it took to solve an MCFP of that size. 
In all cases, the number of communication agents was half the number of task agents, and the weight vector was set to one. 
The results are shown in Table \ref{tab:tab1}. 
We show examples for 30, 35 and 45 agents in Figure \ref{fig:scalability}. 
Furthermore, we can see how the algorithm places network agents in positions that increase the utility function with the number of iterations of the algorithm.

\begin{table}[]
\centering
\begin{tabular}{|l|l|}
\hline
\# Task Agents     & Time to solve MCFP  \\ \hline
A=2  & (0.0005 $\pm$ 0.0001) s      \\ \hline
A=5  & (0.0014 $\pm$ 0.0002) s     \\ \hline
A=10 & (0.0173 $\pm$ 0.0006) s      \\ \hline
A=20 & (0.20	$\pm$ 0.02) s        \\ \hline
A=25 & (0.42 $\pm$ 0.05) s          \\ \hline
A=30 & (1.0 	$\pm$ 0.1) s            \\ \hline
\end{tabular}
\caption{The presented values are the average times calculated for solving a MCFP for different sizes. For each number of task agents we run 50 different scenarios, and solved 10 MCFP for each. In all cases the number of communication agents is half the number of task agents, rounded down on even cases. Experiments were run on an AMD Ryzen Threadripper PRO 3995WX 64-Cores.}
\label{tab:tab1}
\end{table}

\subsection{Mobility}
We lastly show how the algorithm may be used in a dynamic scenario. 
For these simulations, we run two parallel threads. 
The first thread updates the positions of all agents every $\Delta T$ seconds. 
The positions to which task agents are updated follow a trajectory given by a random acceleration model. 
The equations of this model for agent $i \in \mathcal{A}$ are
\begin{align}    
 &\mathbf{v}_i(t+\Delta T) = \mathbf{v}_i(t) +  \mathbf{a}_i(t) \Delta T \\
 &\mathbf{x}_i(t+\Delta T) = \mathbf{x}_i(t) + \mathbf{v}_i(t) \Delta T + \frac{1}{2} \mathbf{a}_i(t) \Delta T^2
\end{align}
where $\mathbf{a}_i(t) \sim \mathcal{N}(0, a I_2)$.
To maintain the density $\rho$ constant, we limit the area in which task agents can navigate. 
In case they hit one of the walls, we change the sign of their velocity to generate a ``bouncing'' effect.
The directions along which positions of network agents are updated are determined by the second thread. 
Letting $\mathbf{d}$ be that direction, the positions are updated as a step in that direction with norm bounded by $v_{max}$.
Here $v_{max}$ is a bound on the speed of the agents. 
The equation for updating the position of agent $i \in \mathcal{I}$ is
\begin{align}
\mathbf{x}_i(t + \Delta T) = \mathbf{x}_i(t) + \min(||\mathbf{d}||,v_{max}) \frac{\mathbf{d}}{||\mathbf{d}||} \Delta T.  
\end{align}
Initial positions for task agents $X_{\mathcal{A}}^0$ are randomly sampled, while initial positions for network agents $X_{\mathcal{I}}^0$ are obtained by finding optimal positions for $X_{\mathcal{A}}^0$ before the simulation starts.
The second thread involves iteratively sampling agent positions for a given instant, solving the MCFP and outputting the directions of local increase for the network agents. 
These directions are the directions $\mathbf{d}$ previously mentioned.

\begin{figure}[h]
    \centering
    \includegraphics[scale=0.25]{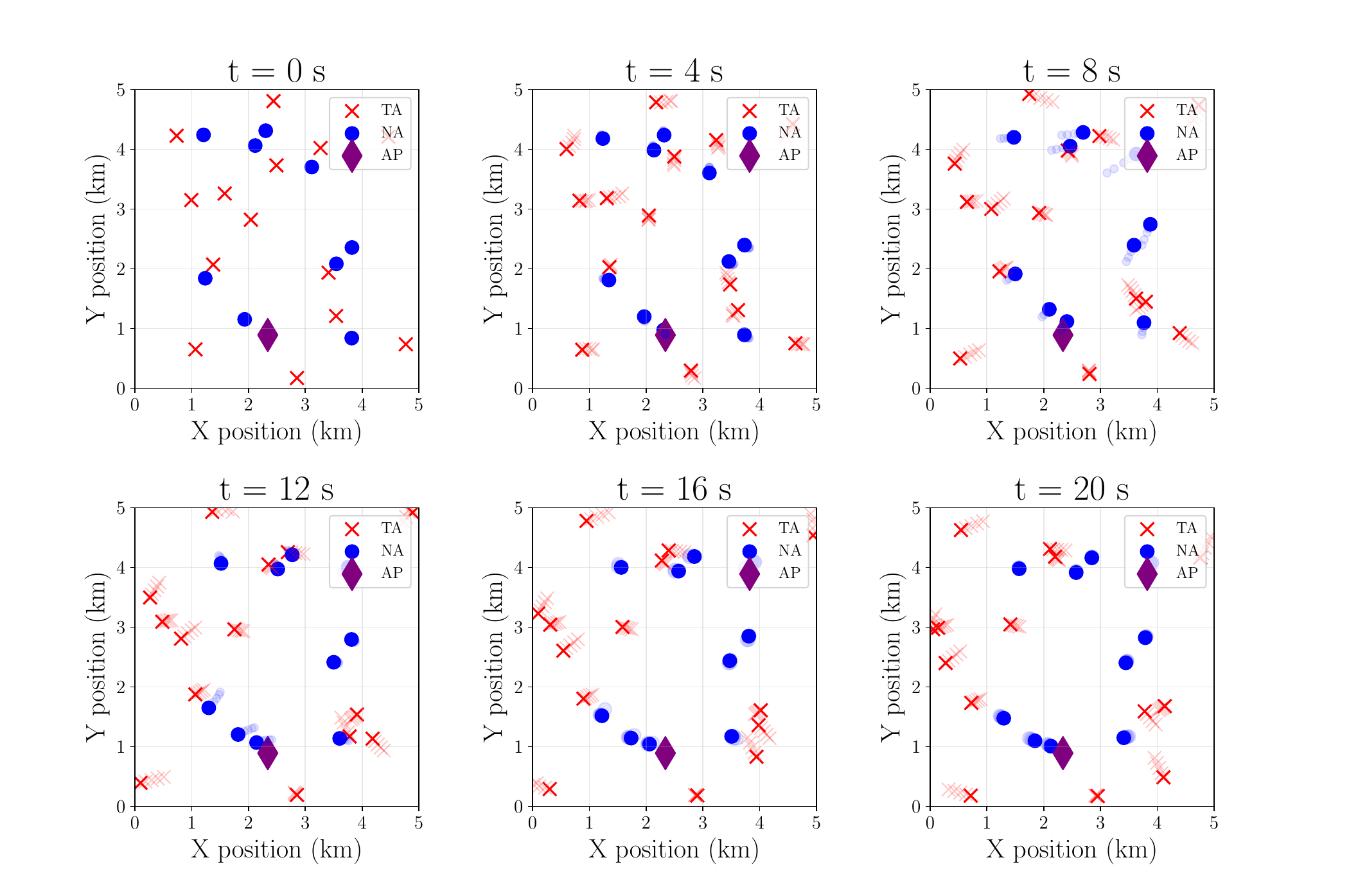}
    \caption{Test for the dynamic scenario simulation. Each frame corresponds to a snap shot of agents configuration for given instants. }
    \label{fig:dynamic}
\end{figure}

We run the scenario simulating case study two. 
The agent that plays the role of the AP is fixed in space. 
We use $\Delta T = 0.2 s$, $a = 0.01 km/s^2$ and $v_{max} = 90 m/s$. 
We simulate for 20 seconds. 
Snapshots of the simulations are shown in Figure \ref{fig:dynamic}. 
We observe how the performance drops in between second 8 and 12 of the simulation. 
That can be attributed to two task agents that separate from the swarm towards the upper right and lower left corner, as can be seen on Figure \ref{fig:dynamic}.
Nonetheless, the network agents close to them move in direction of those distant agents so as to connect them back to the sink, thus recovering the performance function.


\begin{figure}[h]
    \centering
    \includegraphics[scale=0.25]{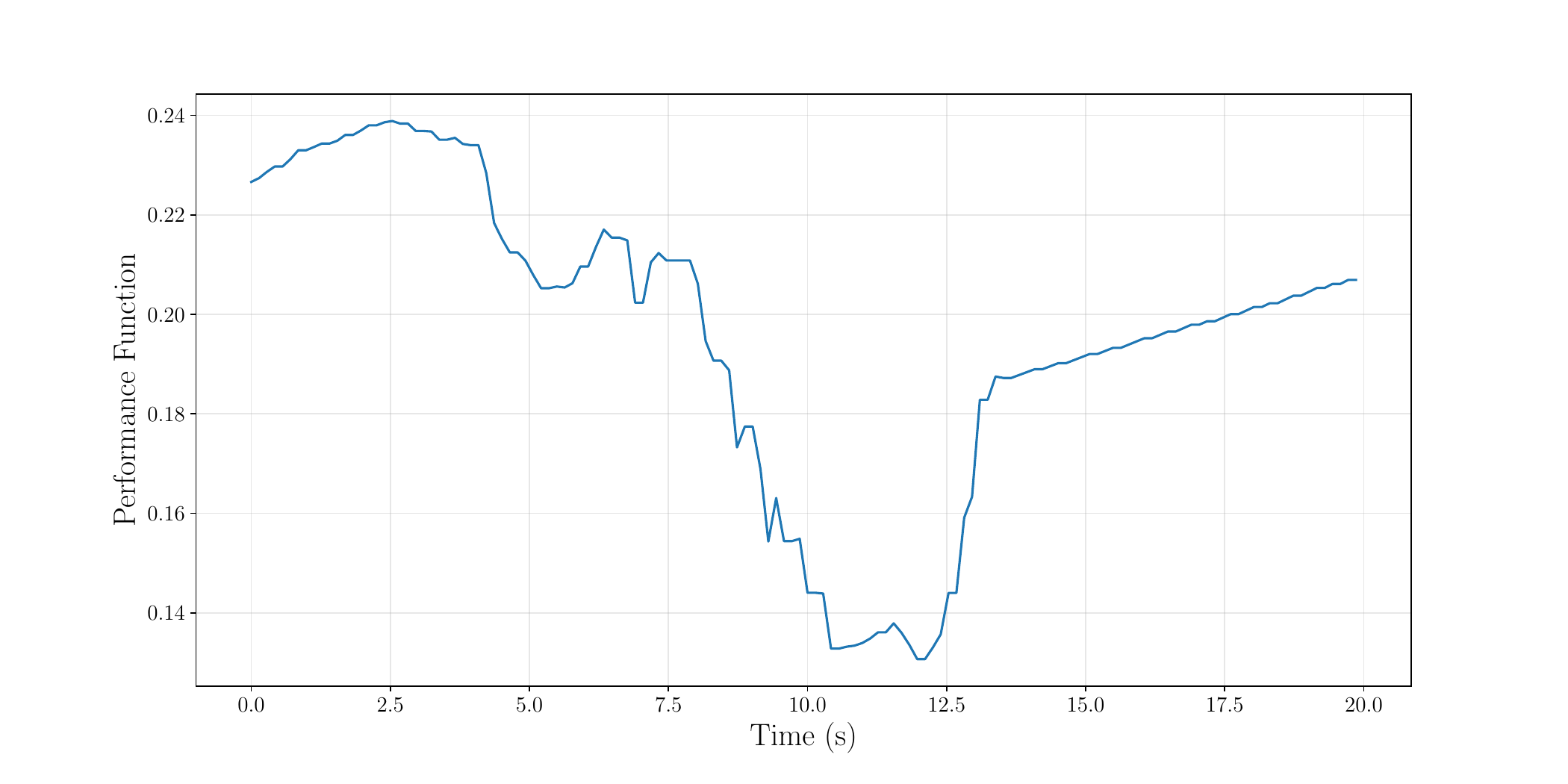}
    \caption{Utility as a function of time for the dynamic scenario simulation.}
    \label{fig:perfDyn}
\end{figure}
\section{Conclusion}

We presented a novel method for tackling the problem of providing dynamic communication infrastructure for a team of agents. 
Our approached leveraged the fact that we can extract directions of local increase of the objective function using the dual solution of a closely related convex optimization problem.
We showed the efficacy of our algorithm on a range of different scenarios in simulation, demonstrating the scalability of the method to large teams of agents. 
Future work will include using learning techniques to decrease the running time of computing the solution to the Multi-Commodity Flow Problem, as well as using different utility functions to model a wider range of practical scenarios.

\bibliographystyle{IEEEtran}
\bibliography{papers}

\end{document}